\documentclass[preprint,12pt]{elsarticle}
\usepackage{amssymb}
\usepackage{subfigure}
\usepackage{amssymb}
\usepackage{amsmath}
\usepackage{amsfonts}
\usepackage{graphicx}
\usepackage{placeins}
\newtheorem{proposition}{Proposition}
\newtheorem{theorem}{Theorem}
\newtheorem{remark}{Remark}
\newtheorem{corollary}{Corollary}
\newtheorem{definition}{Definition}
\journal{Computers \& Fluids}
\begin{document}
\begin{frontmatter}

\title{On a 1/2-equation model of turbulence}
\author[1]{Rui Fang\corref{cor1}}
\ead{ruf10@pitt.edu}
\ead[url]{https://ruf10.github.io}
\cortext[cor1]{Corresponding author: Rui Fang}
\author[2]{Weiwei Han}
\ead{hanweiwei@stu.xjtu.edu.cn}
\author[1]{William Layton}
\ead{wjl@pitt.edu}
\affiliation[1]{organization={Department of Mathematics},
           addressline={University of Pittsburgh}, city={Pittsburgh},
           postcode={15260},
             state={PA},
             country={USA}
             }      
\affiliation[2]{organization={School of Mathematics and
Statistics},
             addressline={ Xi'an Jiaotong University},
             city={Xi'an},
             postcode={710049},
             state={Shaan Xi},
            country={China}
            }
\begin{abstract}
In 1-equation URANS models of turbulence the eddy viscosity is given by
$\nu_{T}=0.55l(x,t)\sqrt{k(x,t)}$ . The length scale $l$ must be pre-specified
and $k(x,t)$ is determined by solving a nonlinear partial differential
equation. We show that in interesting cases the spacial mean of $k(x,t)$
satisfies a simple ordinary differential equation. Using its solution in
$\nu_{T}$ results in a 1/2-equation model. This model has attractive analytic
properties. Further, in comparative tests in 2d and 3d the velocity statistics
produced by the 1/2-equation model are comparable to those of the full
1-equation model.
\end{abstract}
%\item Due to the computational costs of DNS and LES, URANS models are still widely used. This report develops a URANS model  that lowers model complexity and simulation costs preserving current accuracy.

%\item An exact equation for the space mean of the model's TKE is derived. Replacing the TKE with its mean in the eddy viscosity appends a single ODE to the momentum equation. This produces a, so-called, 1/2 equation model.

%\item Model analysis shows that time averaged statistics are well defined, the model of the TKE mean is positive and the model converges to the NSE as the time window decreases. Model analysis also establishes that the resulting eddy viscosity parameterization does not over diffuse the approximate velocity.
%\item The new model's velocity statistics for a 3d flow at statistical equilibrium are comparable to those produced by two full 1-equation models.
%\item Tests of flow statistics for a 2d, higher Re flow not at statistical equilibrium shows the new model's statistics are comparable to both the better 1-equation model and those produced by a finer mesh NSE simulation.


\end{frontmatter}
\section{Introduction}

Unsteady Reynolds averaged Navier Stokes (URANS) models approximate time
averages
\begin{equation}
\overline{u}(x,t):=\frac{1}{\tau}\int_{t-\tau}^{t}u(x,t^{\prime})dt^{\prime
}\text{ with fluctuation }u^{\prime}(x,t):=(u-\overline{u})(x,t)
\label{eq:TimeAverage}%
\end{equation}
of solutions of the Navier-Stokes equations%
\begin{equation}
\nabla\cdot u=0,u_{t}+u\cdot\nabla u-\nu\triangle u+\nabla p=f(x,t),
\label{eqNSE}%
\end{equation}
with the domain, kinematic viscosity, initial and boundary conditions
specified. There are a variety, 0-equation, 1-equation, 2-equation,
more-equation, of useful URANS models with (generally) increasing predictive
ability as model complexity (e.g., number of equations and calibration
parameters) increases. This report studies the extent flow statistics
predicted by 1-equation models can be captured by a 1/2-equation model
(derived in Section 2) which has 0-equation complexity.

The standard URANS approach is to model $\overline{u}(x,t)$ by eddy viscosity%
\begin{equation}
v_{t}+v\cdot\nabla v-\nabla\cdot\left(  \left[  2\nu+\nu_{T}\right]
\nabla^{s}v\right)  +\nabla q=\frac{1}{\tau}\int_{t-\tau}^{t}f(x,t^{\prime
})dt^{\prime}\text{, and }\nabla\cdot v=0.\nonumber
\end{equation}
Here $\nu_{T}=0.55l\sqrt{k}$ is the eddy viscosity. The model representation
of the turbulence length scale\ $l(x,t)$ and turbulent kinetic energy
$k(x,t)\simeq\frac{1}{2}\overline{|u^{\prime}(x,t)|^{2}}$ must be specified.
In Section 2 we show that with kinematic $l=\sqrt{2k}\tau$ the time evolution
of the space-average of $k(x,t)$%
\[
k(t)=\frac{1}{|\Omega|}\int_{\Omega}k(x,t)dx\simeq\frac{1}{|\Omega|}%
\int_{\Omega}\frac{1}{2}\overline{|u^{\prime}(x,t)|^{2}}dx.
\]
can be captured by a single ODE in time%
\begin{equation}
\frac{d}{dt}k(t)+\frac{\sqrt{2}}{2}\tau^{-1}k(t)=\frac{1}{|\Omega|}%
\int_{\Omega}\nu_{T}|\nabla^{s}v|^{2}dx. \label{eq:K(t)Eqn}%
\end{equation}
Using $k(t)$ rather than $k(x,t)$ in $\nu_{T}$ reduces model complexity to
that of a 0-equation model. Section \ref{1/2-equation} proves positivity of
$k(t)$ and boundedness of the 1/2-equation model's kinetic energy and energy
dissipation rate. Proposition \ref{proposition_positivity} in Section
\ref{1/2-equation} shows that when the time window $\tau$\ is sufficiently
small $k(t)\rightarrow0$ (and thus $\nu_{T}$\ also) reducing the model to the
NSE. The other natural limit is whether the model solution converges to a RANS
approximation as $\tau\rightarrow\infty$. Analysis of this question is an open
problem but our preliminary, heuristic analysis suggests it does hold. The
goal of URANS simulations is to give acceptable accuracy at modest cost. One
requirement for this is that the model's eddy viscosity not over dissipate.
This is proven for turbulence in a box in Section 2.4. Section 3 directly
addresses accuracy, comparing 1/2-equation model velocity statistics with
those of 1-equation models. Since simulations of the 1/2-equation model have
reduced complexity compared to 1-equation models, the tests in Section 3
indicate that 1/2-equation model's comparable accuracy makes it worthy of
further study.

\textbf{Related work}. Finite time averaging (\ref{eq:TimeAverage}) is one of
various averages, surveyed by Denaro \cite{D23}, used to develop URANS models.
We select it because it is analytically coherent and computationally feasible.
The equation (\ref{eq:K(t)Eqn}) is derived by space averaging the TKE equation
developed by Prandtl \cite{P45} and Kolmogorov \cite{K42}. The equation for
the spaced average TKE has the simpler form (\ref{eq:K(t)Eqn}) due to the
kinematic turbulence length $l=\sqrt{2k}\tau$, Section 2.2. This $l(x,t)$ was
mentioned by Prandtl, Section 2.1, but developed much later. Our previous work
\cite{KLS21}, \cite{KLS22}, \cite{LM20} has found it to be effective when
boundary layers are not primary and it has been used successfully by\ Teixeira
and Cheinet \cite{TC01}, \cite{TC04} in GFD simulations. Our approach
to\ 1/2-equation models is inspired by the pioneering work of Johnson and King
\cite{JK85}, see also Wilcox \cite{Wilcox} Section 3.7, Johnson \cite{J87}.
This work captured variation of model parameters along a body or channel by
deriving and solving an ODE in $x$, the streamwise direction. We also note
that (\ref{eq:TimeAverage}) means that there is not a sharp separation between
our approach to URANS herein and time filtered large eddy simulation, reviewed
in Pruett \cite{P08}.

\section{The 1/2-equation model}
The 1-equation model is reviewed in Section 2.1 followed by derivation of the
1/2-equation model studied herein in Section 2.2. Analytical properties of the
model are developed in Sections 2.3 and 2.4. Hereafter, we will redefine the
body force to simplify notation, replacing $\frac{1}{\tau}\int_{t-\tau}%
^{t}f(x,t^{\prime})dt^{\prime}$ by $f(x,t)$.

\subsection{Background on 1-equation models}

Averaging the incompressible NSE (\ref{eqNSE}) by (\ref{eq:TimeAverage}) leads
to the exact but non-closed equations for $\overline{u}$:%
\begin{gather*}
\nabla\cdot\overline{u}=0,\overline{u}_{t}+\overline{u}\cdot\nabla\overline
{u}-\nu\triangle\overline{u}+\nabla\overline{p}+\nabla\cdot R(u,u)=f(x,t),\\
\text{where }R(u,u):=\overline{u\otimes u}-\overline{u}\otimes\overline{u}.
\end{gather*}
With few exceptions, URANS models are based on the Boussinesq assumption (that
the action of $R(u,u)$ on $\overline{u}$ is dissipative, \cite{B77}) and the
eddy viscosity hypothesis (that this dissipation can be represented by an
enhanced viscosity $\nu_{T}$, \cite{JLMRZ20}). These yield the model for
$v\simeq\overline{u},$%
\begin{equation}
\nabla\cdot v=0,v_{t}+v\cdot\nabla v-\nu\triangle v-\nabla\cdot\left(  \nu
_{T}\nabla^{s}v\right)  +\nabla q=f(x,t),
\end{equation}
where $q$ is a pressure and $\nabla^{s}v$\ is the symmetric part of the
gradient tensor. Computational experience (now with rigorous mathematical
support \cite{KLS21}, \cite{KLS22}) is that the near wall behavior
$R(u,u)=\mathcal{O}\left(  \left[  \text{wall distance}\right]  ^{2}\right)  $
must be replicated in $\nu_{T}$\ to preclude model over-dissipation. The
turbulent viscosity $\nu_{T}$\ is an expression of the observed increase of
mixing with \textit{\textquotedblleft the intensity of the whirling
agitation\textquotedblright}, \cite{S43}, \cite{B77}, \cite{Darrigol}, p.235.
This results in the dimensionally consistent, Prandtl-Kolmogorov formula%
\begin{gather*}
\nu_{T}=0.55l\sqrt{k}\text{, \ where }l(x,t)=\text{ turbulence length
scale,}\\
k(x,t)\simeq\overline{\frac{1}{2}|u^{\prime}(x,t)|^{2}}\text{ = turbulent
kinetic energy.}%
\end{gather*}
\textit{0-equation models} specify $l$ and relate $k$ back to local changes in
$v(x,t)$. For example, the Baldwin-Lomax \cite{BL78} model uses $k(x,t)\simeq
l^{2}|\nabla\times v(x,t)|^{2}.$ \textit{2-3-...URANS models} solve the
k-equation below for $k(x,t)$. Then they determine $l(x,t)$ indirectly through
the solution of added nonlinear partial differential equations for
dimensionally related turbulent flow statistics.

\textit{1-equation models, }with the notable exception of the Spalart-Alamaras
model, specify $l$ and solve the associated nonlinear PDE for $k(x,t)$
\[
k_{t}+v\cdot\nabla k-\nabla\cdot\left(  \nu_{T}\nabla k\right)  +\frac{1}%
{l}k\sqrt{k}=\nu_{T}|\nabla^{s}v|^{2}.
\]
This is derived by plausible closures of an exact equation for $\overline
{\frac{1}{2}|u^{\prime}|^{2}}$, \cite{CL} p.99, Section 4.4, \cite{D15},
\cite{MP} p.60, Section 5.3 or \cite{Pope} p.369, Section 10.3. Prandtl gave
two descriptions of the physical meaning of $l(x,t)$, e.g. \cite{P45},
\cite{TC01}. The first is that $l(x,t)$ is an average distance turbulent
eddies must go to interact. Walls constrain this distance, leading to
$l=\kappa y$, where $\kappa$ is\ the von Karman constant, Prandtl \cite{P35},
and $y$ is the wall normal distance. The second, kinematic, specification is
the distance a fluctuating turbulent eddy travels in 1 time unit. Their rate
is $|u^{\prime}|\simeq\sqrt{2k(x,t)}$\ leading to a kinematic length scale of%
\begin{equation}
l(x,t)=\sqrt{2k(x,t)}\tau.\tag{Kinematic $l(x,t)$}%
\end{equation}
This is the choice herein and by Kolmogorov \cite{K42} for his 2-equation
model. With $l=\sqrt{2k}\tau$, the 1-equation model becomes $\nu_{T}=\mu
\sqrt{2}k(x,t)\tau$ and
\begin{equation}
\left\{
\begin{array}
[c]{c}%
v_{t}+v\cdot\nabla v-\nabla\cdot\left(  \left[  2\nu+\nu_{T}\right]
\nabla^{s}v\right)  +\nabla q=f(x,t)\text{ and }\nabla\cdot v=0\text{, }\\
k_{t}+v\cdot\nabla k-\nabla\cdot\left(  \nu_{T}\nabla k\right)  +\frac
{\sqrt{2}}{2}\tau^{-1}k=\nu_{T}|\nabla^{s}v|^{2}.
\end{array}
\right.  \label{eq:1EqnModel}%
\end{equation}
At $t=0$ the velocity is initialized by $v(x,0)=v_{0}(x)$. At some $t^{\ast
}\geq0$ the equation for $k(x,t)$ is initialized by $k(x,t^{\ast})=k_{0}(x)$.
We impose standard no-slip boundary conditions\footnote{When $\nu
_{T}=\mathcal{O}\left(  \left[  \text{wall distance}\right]  ^{2}\right)  $
there is a serious analytical question about the meaning of traces of $k(x,t)$
on the domain boundary. This question does not arise for the 1/2-equation
model.} at walls: $v(x,t)=0,k(x,t)=0$ on $\partial\Omega.$

\subsection{Derivation of the 1/2-equation model}

We select the kinematic length scale $l=\sqrt{2k}\tau$ yielding $\nu_{T}%
=\sqrt{2}\mu k(x,t)\tau$ and (\ref{eq:1EqnModel}). The 1/2-equation model
begins with the space average of the k-equation in (\ref{eq:1EqnModel}). Let%
\[
k(t):=\frac{1}{|\Omega|}\int_{\Omega}k(x,t)dx\text{ and }\varepsilon
(t):=\frac{1}{|\Omega|}\int_{\Omega}\nu_{T}|\nabla^{s}v|^{2}dx.
\]
Averaging the k-equation over $\Omega$ gives%

\[
\frac{d}{dt}\frac{1}{|\Omega|}\int_{\Omega}kdx+\frac{1}{|\Omega|}\int_{\Omega
}\nabla\cdot(vk)-\nabla\cdot\left(  \nu_{T}\nabla k\right)  dx+\frac{\sqrt{2}%
}{2}\tau^{-1}\frac{1}{|\Omega|}\int_{\Omega}kdx=\varepsilon(t)
\]
The second and third terms vanish since $k(x,t)=0$ on $\partial\Omega:$
\begin{equation}
\left\{
\begin{array}
[c]{c}%
\int_{\Omega}\nabla\cdot(vk)dx=\int_{\partial\Omega}(v\cdot n)kd\sigma=0,\\
\int_{\Omega}\nabla\cdot\left(  \nu_{T}\nabla k\right)  dx=\int_{\partial
\Omega}\sqrt{2}\mu k(x,t)\tau\nabla k\cdot nd\sigma=0.
\end{array}
\right.  \label{eq:BoundaryIntegrals}%
\end{equation}
In more detail, both integrals are zero for internal flows ($v=0$ and $k=0$ on
$\partial\Omega$), shear flows ($v\cdot n=0$ and $k=0$\ on $\partial\Omega$)
and under periodic boundary conditions. For these 3 cases, $k(t)$ satisfies
exactly%
\begin{equation}
\frac{d}{dt}k(t)+\frac{\sqrt{2}}{2}\tau^{-1}k(t)=\varepsilon(t)\text{ with
}k(t^{\ast})\neq0\text{\ given.}%
\end{equation}
One further\ model refinement is needed near walls. With $k=k(t)$, $\nu
_{T}=\sqrt{2}\mu k(t)\tau$ does not vanish at walls. Recall $y=$\ \textit{wall
normal distance }and\textit{ }$\kappa=$ \textit{von Karman constant}. Since
$\nu_{T}$\ should replicate the $\mathcal{O}(y^{2})$ near wall asymptotics of
$R(u,u)$, we adjust $\nu_{T}$\ at walls with a multiplier $\left(  \kappa
y/L\right)  ^{2}$. We thus have the 1/2-equation model%
\begin{gather}
\nu_{T}=\sqrt{2}\mu k(t)\tau\text{ for periodic boundaries,}\nonumber\\
\nu_{T}=\sqrt{2}\mu k(t)\left(  \kappa\frac{y}{L}\right)  ^{2}\tau\text{ for
no-slip and shear boundaries,}\nonumber\\
v_{t}+v\cdot\nabla v-\nabla\cdot\left(  \left[  2\nu+\nu_{T}\right]
\nabla^{s}v\right)  +\nabla q=f(x,t)\text{ and }\nabla\cdot v=0\text{,
}\nonumber\\
\frac{d}{dt}k(t)+\frac{\sqrt{2}}{2}\tau^{-1}k(t)=\varepsilon(t)\text{
}.\label{eqHalfEqnModel}%
\end{gather}

\begin{remark}
\textbf{Channel flows} have walls but also inflow and outflow boundaries.
Outflow BCs are non-physical and selected to do minimal harm to the upstream
approximation. Inflow values for $v(x,t),k(x,t)$ are needed and must be
specified from measurements%
\[
v\cdot n=v_{IN}\text{ , }v\cdot\tau=0,\text{ }k=k_{IN}\text{ on inflow
boundary }\Gamma_{IN}\text{.}%
\]
With these known the volume averaged inflow part of the convection integral is
calculable:%
\[
\frac{1}{|\Omega|}\int_{\partial\Omega}(v\cdot n)kd\sigma=\frac{1}{|\Omega
|}\int_{\Gamma_{IN}}v_{IN}\text{ }k_{IN}\text{ }d\sigma\text{ }+\text{ outflow
integral}.
\]
The volume averaged diffusion integral becomes $\int_{\Gamma_{IN}}%
k_{IN}(x)\nabla k\cdot nd\sigma$ and is not exactly calculable due to $\nabla
k$ being unknown on $\Gamma_{IN}$. When diffusion of TKE across the inflow
boundary is much smaller than convection (a plausible but untested
hypothesis), the second integral is negligible. Under these conditions the
1/2-equation model for channel flow is%
\[
\frac{d}{dt}k(t)+\frac{\sqrt{2}}{2}\tau^{-1}k(t)=\varepsilon(t)-\frac
{1}{|\Omega|}\int_{\Gamma_{IN}}u_{IN}\text{ }k_{IN}\text{ }d\sigma\text{
}-\text{ outflow integral.}%
\]

\end{remark}

\begin{remark}
[URANS not RANS]The approach to replacing $k(x,t)$ with $k(t)$ means that the
model so derived is essentially a URANS model not a RANS model. Nevertheless,
one can ask what happens for flows at statistical equilibrium. At statistical
equilibrium, (and under periodic BCs) $k(t)$ approximately satisfies%
\[
\frac{\sqrt{2}}{2}\tau^{-1}k(t)\simeq\frac{1}{|\Omega|}\int_{\Omega}\sqrt
{2}\mu k(t)\tau|\nabla^{s}v|^{2}dx
\]
and $k(t)$ can be approximately cancelled from the equation. However, implicit
$k$ dependence remains since $v=v(x,t;k)$. At statistically equilibrium the
value of $k$ is thus (approximately) determined by%
\[%
\begin{array}
[c]{cc}%
\text{solve for }k\text{:} & \frac{1}{|\Omega|}\int_{\Omega}|\nabla
^{s}v(x,t;k)|^{2}dx\simeq\frac{1}{2\mu}\tau^{-2}\\
\text{subject to:} & \left\{
\begin{array}
[c]{c}%
v_{t}+v\cdot\nabla v-\nabla\cdot\left(  \left[  2\nu+\sqrt{2}\mu k\tau\right]
\nabla^{s}v\right)  +\nabla q=f(x)\\
\nabla\cdot v=0
\end{array}
\right.
\end{array}
\]

\end{remark}

\subsection{1/2-equation Model: Basic Estimates}

\label{1/2-equation}

We assume here that the model (\ref{eqHalfEqnModel}) under periodic or no slip
boundary conditions has a solution which is smooth enough for standard energy
estimates. While proving this for the continuum model is an open problem,
existence certainly holds for its FEM discretization. We first establish
positivity of $k(t)$ and model convergence to the NSE for small $\tau$.
Analysis of model behavior for $\tau\rightarrow\infty$\ is an open problem.

\begin{proposition}
\label{proposition_positivity} Consider the model (\ref{eqHalfEqnModel}) under
no-slip or periodic boundary conditions. The following hold.

\textbf{Positivity:} If $k(t^{\ast})>0$ then $k(t)>0$ for all $t\geq t^{\ast}$.

\textbf{Model convergence to NSE:} There is a $\tau_{0}=\tau_{0}(data)>0$ such
that for $\tau<\tau_{0}$, $k(t)\rightarrow0$ and $\nu_{T}\rightarrow0$
exponentially in $t$.
\end{proposition}

\textbf{proof:} \textit{The k-equation can be rewritten as}%
\[
\frac{d}{dt}k(t)+\left[  \frac{\sqrt{2}}{2}\tau^{-1}-\tau\frac{\mu}{L^{2}%
}\frac{1}{|\Omega|}\int_{\Omega}y^{2}|\nabla^{s}v(x,t)|^{2}dx\right]
k(t)=0\text{.}%
\]
\textit{Let}%
\[
A(t)=\int_{0}^{t}\left[  \frac{\sqrt{2}}{2}\tau^{-1}-\tau\frac{\mu}{L^{2}%
}\frac{1}{|\Omega|}\int_{\Omega}y^{2}|\nabla^{s}v(x,t^{\prime})|^{2}dx\right]
dt^{\prime}.
\]
\textit{The standard velocity energy estimate implies }$\int y^{2}|\nabla
^{s}v(x,t)|^{2}dx\in L^{1}(0,T)$\textit{. Thus so is the term in brackets and
its antiderivative }$\mathit{A(t)}$\textit{ is well defined. The solution is
then }$k(t)=\exp\{-[A(t)-A(t^{\ast})]\}k(t^{\ast})$\textit{. Thus positivity
of }$k(t)$\textit{ follows. When }$\mathit{A(t)<0}$\textit{, which occurs for
}$\tau<\tau_{0}$\textit{, }$\mathit{k(t)}$\textit{ decays as claimed,
completing the proof. }

Next we establish an energy equality suggested by the NSE's kinetic energy
balance, rewritten in terms of $\int\frac{1}{2}|\overline{u}|^{2}+\frac{1}%
{2}|u^{\prime}|^{2}dx$, and associated \'{a} priori bounds on kinetic energy
and energy dissipation.

\begin{proposition}
Consider the model (\ref{eqHalfEqnModel}) under no-slip or periodic boundary
conditions. Sufficiently regular model solutions satisfy the energy equality
\begin{gather}
\frac{d}{dt}\left[  \frac{1}{|\Omega|}\int_{\Omega}\frac{1}{2}|v(x,t)|^{2}%
dx+k(t)\right]  +\label{EnergyEst}\\
\left[  \frac{1}{|\Omega|}\int_{\Omega}\nu|\nabla^{s}v(x,t)|^{2}dx+\frac
{\sqrt{2}}{2}\tau^{-1}k(t)\right]  =\frac{1}{|\Omega|}\int_{\Omega}f(x,t)\cdot
v(x,t)dx.\nonumber
\end{gather}
Suppose $k(t^{\ast})>0$. With $C=C(data)=C(f,u(x,0),k(t^{\ast}),\nu
,\tau)<\infty,$ the following uniform in $T$\ bounds on energy and dissipation
rates hold%
\begin{align*}
\frac{1}{|\Omega|}\int_{\Omega}\frac{1}{2}|v(x,T)|^{2}dx  &  \leq C,\\
\frac{1}{T}\int_{0}^{T}\left\{  \frac{1}{|\Omega|}\int_{\Omega}[\nu+\nu
_{T}]|\nabla^{s}v(x,t)|^{2}dx\right\}  dt  &  \leq C,\\
\frac{1}{|\Omega|}\int_{\Omega}\frac{1}{2}|v(x,T)|^{2}dx+k(T)  &  \leq C\text{
},\\
\frac{1}{T}\int_{0}^{T}\left\{  \frac{1}{|\Omega|}\int_{\Omega}\nu|\nabla
^{s}v(x,t)|^{2}dx+\frac{\sqrt{2}}{2}\tau^{-1}k(t)\right\}  dt  &  \leq C\text{
}.
\end{align*}

\end{proposition}

\textbf{proof:}\textit{ Take the inner product of the momentum equation with
}$v(x,t)$\textit{, apply the divergence theorem. This gives }%
\begin{gather*}
\frac{d}{dt}\frac{1}{|\Omega|}\int_{\Omega}\frac{1}{2}|v(x,t)|^{2}dx+\\
\frac{1}{|\Omega|}\int_{\Omega}[\nu+\nu_{T}]|\nabla^{s}v(x,t)|^{2}dx=\frac
{1}{|\Omega|}\int_{\Omega}f(x,t)\cdot v(x,t)dx.
\end{gather*}
\textit{Since }$k(t)\geq0,$ $\nu_{T}|\nabla^{s}v(x,t)|^{2}\geq0$ \textit{and
the }$\nu_{T}$\textit{\ term can be dropped (for the kinetic energy bound)
then reinserted (for the dissipation bound). Differential inequalities imply
that, uniformly in} $T$\textit{,}
%\[
%\frac{1}{|\Omega|}\int_{\Omega}\frac{1}{2}|v(x,T)|^{2}dx\leq C(data)%%<\infty
%,
%\frac{1}{T}\int_{0}^{T}\frac{1}{|\Omega|}\int_{\Omega}[\nu+\nu_{T}%
%]|\nabla^{s}v(x,t)|^{2}dxdt\leq C(data)<\infty.\text{ }%
%\]
%Rui changed it to two lines, since its overfull:%
\begin{align*}
\frac{1}{|\Omega|}\int_{\Omega}\frac{1}{2}|v(x,T)|^{2}dx &  \leq
C(data)<\infty,\\
\frac{1}{T}\int_{0}^{T}\frac{1}{|\Omega|}\int_{\Omega}[\nu+\nu_{T}]|\nabla
^{s}v(x,t)|^{2}dxdt &  \leq C(data)<\infty.
\end{align*}
\textit{Adding the k-equation gives}%
\begin{gather}
\frac{d}{dt}\left[  \frac{1}{|\Omega|}\int_{\Omega}\frac{1}{2}|v(x,t)|^{2}%
dx+k(t)\right]  +\\
\left[  \frac{1}{|\Omega|}\int_{\Omega}\nu|\nabla^{s}v(x,t)|^{2}dx+\frac
{\sqrt{2}}{2}\tau^{-1}k(t)\right]  =\frac{1}{|\Omega|}\int_{\Omega}f(x,t)\cdot
v(x,t)dx.\nonumber
\end{gather}
\textit{As above, standard differential inequalities again imply that,
uniformly in }$\mathit{T}$%
\begin{align*}
\frac{1}{|\Omega|}\int_{\Omega}\frac{1}{2}|v(x,T)|^{2}dx+k(T) &  \leq C,\\
\frac{1}{T}\int_{0}^{T}\left\{  \frac{1}{|\Omega|}\int_{\Omega}\nu|\nabla
^{s}v(x,t)|^{2}dx+\frac{\sqrt{2}}{2}\tau^{-1}k(t)\right\}  dt &  \leq C,
\end{align*}
\textit{completing the proof}.

One consequence is the following result on time averaged equilibrium of the k-equation.

\begin{corollary}
As $T\rightarrow\infty$\ there holds%
\[
\frac{1}{T}\int_{0}^{T}\frac{1}{|\Omega|}\int_{\Omega}\nu_{T}dxdt=\frac
{\mu\tau}{T}\int_{0}^{T}k(t)dt=\frac{\sqrt{2}\mu\tau^{2}}{T}\int_{0}%
^{T}\varepsilon(t)dt+\mathcal{O}\left(  \frac{1}{T}\right)  .
\]

\end{corollary}

\textbf{proof:} \textit{Time averaging the }$k-$\textit{equation gives}%
\[
\frac{1}{T}\left(  k(T)-k(0)\right)  +\frac{\sqrt{2}}{2}\tau^{-1}\frac{1}%
{T}\int_{0}^{T}k(t)dt=\frac{1}{T}\int_{0}^{T}\varepsilon(t)dxdt\text{.}%
\]
\textit{The first term is }$O\left(  1/T\right)  $\textit{ due to the a priori
bounds. Rearranging this gives the claimed result}%
\[
\mu\tau\frac{1}{T}\int_{0}^{T}k(t)dt=\sqrt{2}\mu\tau^{2}\frac{1}{T}\int
_{0}^{T}\varepsilon(t)dt+\mathcal{O}\left(  \frac{1}{T}\right)  \text{.}%
\]

\subsection{Energy dissipation rate: turbulence in a box}

We show next that the model does not over dissipate body forced flow with
periodic boundary conditions, often called turbulence in a box. These
estimates use the \'{a} priori bounds in Section 2.3 but require a small
amount of extra notation.

\begin{definition}
The model energy dissipation rate is%
\[
\varepsilon_{\text{model}}(t):=\frac{1}{|\Omega|}\int_{\Omega}\nu|\nabla
^{s}v(x,t)|^{2}dx+\frac{\sqrt{2}}{2}\tau^{-1}k(t).
\]
The scale of the body force $F$, large velocity scale $U$, length scale $L$
and large scale turnover time $T^{\ast}$ are, respectively,
\begin{align*}
F &  :=\sqrt{\frac{1}{|\Omega|}\int_{\Omega}|f(x)|^{2}dx},\\
U &  :=\sqrt{\lim\sup_{T\rightarrow\infty}\frac{1}{T}\int_{0}^{T}\left[
\frac{1}{|\Omega|}\int_{\Omega}|v(x,t)|^{2}dx\right]  dt},\\
L &  :=\min \Bigl\{  |\Omega|^{1/3},\frac{F}{\sup_{x\in\Omega}|\nabla^{s}%
f(x)|},\frac{F}{\sqrt{\frac{1}{|\Omega|}\int_{\Omega}|\nabla^{s}f(x)|^{2}dx}%
},\\
&\sqrt{\frac{F}{\sqrt{\frac{1}{|\Omega|}\int_{\Omega}|\triangle f(x)|^{2}dx}%
}}\Bigl\}  ,\\
T^{\ast} &  :=\frac{L}{U}.
\end{align*}

\end{definition}

The large velocity scale $U$ and length scale $L$ are well defined due to
Proposition 2.3.

\begin{theorem}
Consider the model (\ref{eqHalfEqnModel}) in $3d$ subject to periodic boundary
conditions and with $\nabla\cdot f(x)=0$. The time averaged energy dissipation
rate of the model satisfies%
\[
\left[  1-\sqrt{2}\mu\left(  \frac{\tau}{T^{\ast}}\right)  ^{2}\right]
\lim\sup_{T\rightarrow\infty}\frac{1}{T}\int_{0}^{T}\varepsilon_{\text{model}%
}(t)dt\leq2(1+\mathcal{R}e^{-1})\frac{U^{3}}{L}.
\]
With $\mu=0.55,$ if $\frac{\tau}{T^{\ast}}\leq\allowbreak0.8$ then%
\[
\lim\sup_{T\rightarrow\infty}\frac{1}{T}\int_{0}^{T}\varepsilon_{\text{model}%
}(t)dt\leq4(1+\mathcal{R}e^{-1})\frac{U^{3}}{L}.
\]

\end{theorem}

\textbf{proof:} \textit{\ Let }$\phi(T)$\textit{\ denote a generic, bounded,
positive function with }$\phi(T)\rightarrow0$\textit{\ as }$T\rightarrow
\infty$\textit{. Consider the energy estimate (\ref{EnergyEst}) above (which
establishes that }$\varepsilon_{\text{model}}(t)$ is defined correctly).
\textit{Time averaging (\ref{EnergyEst}) gives}%
\begin{gather}
\frac{1}{T}\left[  \left(  \frac{1}{|\Omega|}\int_{\Omega}\frac{1}%
{2}|v(x,T)|^{2}dx+k(T)\right)  -\left(  \frac{1}{|\Omega|}\int_{\Omega}%
\frac{1}{2}|v(x,0)|^{2}dx+k(0)\right)  \right]  +\\
\frac{1}{T}\int_{0}^{T}\varepsilon_{\text{model}}(t)dt=\frac{1}{T}\int_{0}%
^{T}\left[  \frac{1}{|\Omega|}\int_{\Omega}f(x)\cdot v(x,t)dx\right]
dt.\nonumber
\end{gather}
\textit{From the \'{a} priori bounds in Proposition 2.3, the first term is
}$\mathcal{O}\left(  1/T\right)  $\textit{; the second term is the time
average of }$\varepsilon_{\text{model}}$\textit{. We thus have}%
\begin{align}
\frac{1}{T}\int_{0}^{T}\varepsilon_{\text{model}}(t)dt &  =\mathcal{O}\left(
\frac{1}{T}\right)  +\frac{1}{T}\int_{0}^{T}\left[  \frac{1}{|\Omega|}%
\int_{\Omega}f(x)\cdot v(x,t)dx\right]  dt\nonumber\\
&  \leq\mathcal{O}\left(  \frac{1}{T}\right)  +UF+\phi(T).\label{EqFirstEst}%
\end{align}
\textit{Next take the inner product of (\ref{eqHalfEqnModel}) with }%
$f(x)$\textit{, integrate by parts, use }$\nabla\cdot f=0$\textit{ and time
average over }$0\leq t\leq T$\textit{. This gives}%
\begin{gather}
F^{2}=\frac{1}{T}\int_{0}^{T}\left(  \frac{1}{|\Omega|}\int_{\Omega
}[v(x,T)-v(x,0)]\cdot f(x)dx\right)  dt\label{F2Equation}\\
-\frac{1}{T}\int_{0}^{T}\left(  \frac{1}{|\Omega|}\int_{\Omega}v(x,t)\otimes
v(x,t):\nabla f(x)dx\right)  dt+\nonumber\\
\frac{1}{T}\int_{0}^{T}\left[  \frac{1}{|\Omega|}\int_{\Omega}[\nu+\nu
_{T}]\nabla^{s}v(x,t):\nabla^{s}f(x)dx\right]  dt.\nonumber
\end{gather}
\textit{The first two terms in the RHS are shared with the NSE. Standard
estimates for those terms shared with the NSE from the pioneering work of
Doering, Foias and Constantin \cite{DF02}, \cite{DC92} are }$\mathcal{O}%
\left(  \frac{1}{T}\right)  +\frac{F}{L}U^{2}+\phi(T)$\textit{. Consider the
last term on the RHS. Application of the space-time Cauchy-Schwarz-Young
inequality to it gives, for any }$0<\beta<1$\textit{,}%
\begin{gather}
\frac{1}{T}\int_{0}^{T}\left[  \frac{1}{|\Omega|}\int_{\Omega}[\nu+\nu
_{T}]\nabla^{s}v(x,t):\nabla^{s}f(x)dx\right]  dt\leq\label{eqProofStep}\\
\leq\frac{1}{T}\int_{0}^{T}\frac{1}{|\Omega|}\int_{\Omega}-\nu v:\triangle
fdxdt+\frac{1}{T}\int_{0}^{T}\frac{1}{|\Omega|}\int_{\Omega}\nu_{T}\nabla
^{s}v:\nabla^{s}fdxdt\nonumber\\
\leq\frac{\nu UF}{L^{2}}+\frac{\beta}{2}\frac{F}{UT}\int_{0}^{T}%
\varepsilon(t)dt+\frac{1}{2\beta}\frac{UF}{L^{2}T}\int_{0}^{T}\frac{1}%
{|\Omega|}\int_{\Omega}\mu k(t)\tau dxdt+\phi(T).\nonumber
\end{gather}
\textit{Inserting the estimate of the time average of }$k(t)$ \textit{from
Corollary 2.1 in (\ref{eqProofStep}) gives}%
\begin{gather*}
\frac{1}{T}\int_{0}^{T}\left[  \frac{1}{|\Omega|}\int_{\Omega}[\nu+\nu
_{T}]\nabla^{s}v(x,t):\nabla^{s}f(x)dx\right]  dt\leq\nu U\frac{F}{L^{2}}+\\
+\frac{\beta}{2}\frac{F}{U}\frac{1}{T}\int_{0}^{T}\varepsilon_{\text{model}%
}(t)dt+\frac{\sqrt{2}\mu\tau^{2}}{2\beta}U\frac{F}{L^{2}}\frac{1}{T}\int
_{0}^{T}\varepsilon(t)dt+\mathcal{O}\left(  \frac{1}{T}\right)  +\phi(T).
\end{gather*}
\textit{Use this and the previously derived, }$\mathcal{O}\left(  \frac{1}%
{T}\right)  +\frac{F}{L}U^{2}+\phi(T)$\textit{, estimates for the RHS terms in
the equation for }$F^{2}$\textit{ (\ref{F2Equation}). This gives, for any
}$0<\beta<1$\textit{,}%
\begin{align*}
F^{2} &  \leq\nu U\frac{F}{L^{2}}+\frac{F}{L}U^{2}+\frac{\beta}{2}\frac{F}%
{U}\frac{1}{T}\int_{0}^{T}\varepsilon_{\text{model}}(t)dt\\
&  +\frac{\sqrt{2}\mu\tau^{2}}{2\beta}U\frac{F}{L^{2}}\frac{1}{T}\int_{0}%
^{T}\varepsilon(t)dt+\mathcal{O}\left(  \frac{1}{T}\right)  +\phi(T).
\end{align*}
\textit{Therefore, we have the key inequality in estimating }$\frac{1}{T}%
\int_{0}^{T}\varepsilon_{\text{model}}(t)dt$\textit{: }%
\begin{align*}
UF &  \leq\nu\frac{U^{2}}{L^{2}}+\frac{U^{3}}{L}+\frac{\beta}{2}\frac{1}%
{T}\int_{0}^{T}\varepsilon_{\text{model}}(t)dt\\
&  +\frac{\sqrt{2}\mu\tau^{2}}{2\beta}\frac{U^{2}}{L^{2}}\frac{1}{T}\int
_{0}^{T}\varepsilon(t)dt+\mathcal{O}\left(  \frac{1}{T}\right)  +\phi(T).
\end{align*}
\textit{The first term on the RHS }$\nu\frac{U^{2}}{L^{2}}$\textit{ is
rewritten as }$\nu\frac{U^{2}}{L^{2}}=\frac{\nu}{LU}\frac{U^{3}}{L}%
=Re^{-1}\frac{U^{3}}{L}.$\textit{ \ Insert this on the RHS of the above and
replace }$\mathit{UF}$\textit{ in (\ref{EqFirstEst}) by the last bound. This
yields}%
\begin{align*}
\frac{1}{T}\int_{0}^{T}\varepsilon_{\text{model}}(t)dt &  \leq\frac{U^{3}}%
{L}+\mathcal{R}e^{-1}\frac{U^{3}}{L}+\frac{\beta}{2}\frac{1}{T}\int_{0}%
^{T}\varepsilon_{\text{model}}(t)dt+\\
&  +\frac{\sqrt{2}\mu}{2\beta}\left\{  \frac{\tau^{2}U^{2}}{L^{2}}\right\}
\frac{1}{T}\int_{0}^{T}\varepsilon(t)dt+\mathcal{O}\left(  \frac{1}{T}\right)
+\phi(T).
\end{align*}
\textit{The multiplier in braces }$\frac{\tau^{2}U^{2}}{L^{2}}=\left(
\frac{\tau}{T^{\ast}}\right)  ^{2}.$ \textit{Pick }$\beta=1$\textit{ and use}%
\[
\frac{1}{T}\int_{0}^{T}\varepsilon(t)dt=\frac{1}{T}\int_{0}^{T}\frac{\sqrt{2}%
}{2}\tau^{-1}k(t)dt+{\small O}\left(  1\right)  \leq\frac{1}{T}\int_{0}%
^{T}\varepsilon_{\text{model}}(t)dt+\phi(T).
\]
\textit{ Proposition 2.3 shows that }$\frac{1}{T}\int_{0}^{T}\varepsilon
_{\text{model}}(t)dt$\textit{ is bounded uniformly in }$T$\textit{. Thus its
limit superior as }$T\rightarrow\infty$ exists. \textit{This (plus an
arithmetic calculation) completes the proof:}%
\[
\left[  1-\sqrt{2}\mu\left(  \frac{\tau}{T^{\ast}}\right)  ^{2}\right]
\lim\sup_{T\rightarrow\infty}\frac{1}{T}\int_{0}^{T}\varepsilon_{\text{model}%
}(t)dt\leq2(1+\mathcal{R}e^{-1})\frac{U^{3}}{L}.
\]

\section{Testing the 1/2-equation model}

We test how close the 1/2-equation model velocity statistics are to volume
averaged velocity statistics produced by the 1-equation model with Prandtl's
classical $l=0.41y$ and with the kinematic turbulence length scale
$l=\sqrt{2k(x,t)}\tau$. Since the 1/2-equation model's $k(t)$ allows temporal
variations, our intuition is that a time independent body force, leading to a
flow where statistical equilibrium is expected, is a non-trivial test.

\subsection{Flow statistics}

Evaluation of flow statistics means comparing plots of $1d$ curves of
aggregate velocity-based quantities. We calculate the time evolution of the
four velocity statistics: the Taylor microscale (an average velocity length
scale), kinetic energy of the mean flow, enstrophy (aggregate vorticity) and
the model approximation to the turbulent intensity:%
\[%
\begin{array}
[c]{ccc}%
\text{Taylor Microscale} & : & \lambda_{Taylor}:=\frac{1}{15}\left(
\frac{\frac{1}{|\Omega|}\int_{\Omega}|\nabla^{s}v(x,t)|^{2}dx}{\frac
{1}{|\Omega|}\int_{\Omega}|v(x,t)|^{2}dx}\right)  ^{-1/2}\\
\text{Kinetic Energy/Volume} & : & E(t):=\frac{1}{|\Omega|}\int_{\Omega}%
\frac{1}{2}|v(x,t)|^{2}dx\\
\text{Enstrophy/Volume} & : & Ens(t):=\frac{1}{|\Omega|}\int_{\Omega}\frac
{1}{2}|\nabla\times v(x,t)|^{2}dx\\
\text{Turbulence Intensity} & : & I_{\text{model}}(t)=\frac{\frac{1}{|\Omega
|}\int_{\Omega}2k(x,t)dx}{\frac{1}{|\Omega|}\int_{\Omega}2k(x,t)+|v(x,t)|^{2}%
dx}%
\end{array}
\]

\subsection{The 3d test problem}

We examined the classical Taylor-Couette flow between counter-rotating
cylinders with no-slip BCs. We used FEniCSx with the computational environment
DOLFINx/0.5.2. We compared three models: the 1/2-equation model in Section
\ref{1/2-equation}, the 1-equation model in equation (\ref{eq:1EqnModel}), and
Prandtl's classical model with $l=0.41y$, where $y$ is the wall normal
distance. We do not know if the von Karman constant is the correct calibration
in the multiplier $(\kappa y/L)^{2}$ . In the 3d test we tested $\kappa=1$ in
the multiplier $(\kappa y/L)^{2}=(y/L)^{2}$. We used the backward Euler time
discretization for both the momentum and k-equation plus a time filter from
\cite{GL18} for the momentum equation to increase time accuracy and
anti-diffuse the implicit method. We used the Taylor-Hood ($P2-P1$) element
pair for the momentum equation in all cases. For the 1-equation model
simulations we used $P1$ Lagrange elements for the k-equation. The
unstructured mesh was generated with GMSH, with GMSH target mesh size
parameter $lc$ $=0.04$. The domain is given by
\[
\Omega=\{(x,y,z):r_{inner}^{2}\leq x^{2}+y^{2}\leq r_{outer}^{2},0\leq z\leq
z_{max}\}
\]
with $r_{inner}=0.833$ and $r_{outer}=1$ and $z_{max}=2.2$. Periodic boundary
conditions were imposed in the $z$ direction. The outer cylinder remained
stationary, while the inner cylinder rotation drove the flow. The angular
velocity of the inner cylinder, $\omega_{inner}$, started at 0 at $T=0$, and
gradually increased until fully turned on with $\omega_{inner}=9$ at time
$T=5$. We chose the final time $T=30$. The time scale was set to be $\tau
=0.1$, and timestep $\Delta t=10^{-2}$. We set $\nu=10^{-3}$, took $U=$ inner
cylinder speed and $L=1.0-0.833=0.167$, the cylinder gap, yielding
$\mathcal{R}e=1.5\times10^{3}$. The radius ratio $\eta$ and Taylor number $Ta$
are
\[
\eta:=\frac{r_{inner}}{r_{outer}}=0.833\text{ and }Ta:=4\mathcal{R}e^{2}%
\frac{1-\eta}{1+\eta}\simeq8\times10^{5}.
\]
Figure 1, p. 156 in \cite{ALS86} (see also \cite{GLS2016}) indicates the
physical flow is expected to have turbulent Taylor cells for these parameters.
\FloatBarrier\begin{figure}[ptbh]
\begin{minipage}{0.5\textwidth}
\centering
\includegraphics[width=\linewidth]{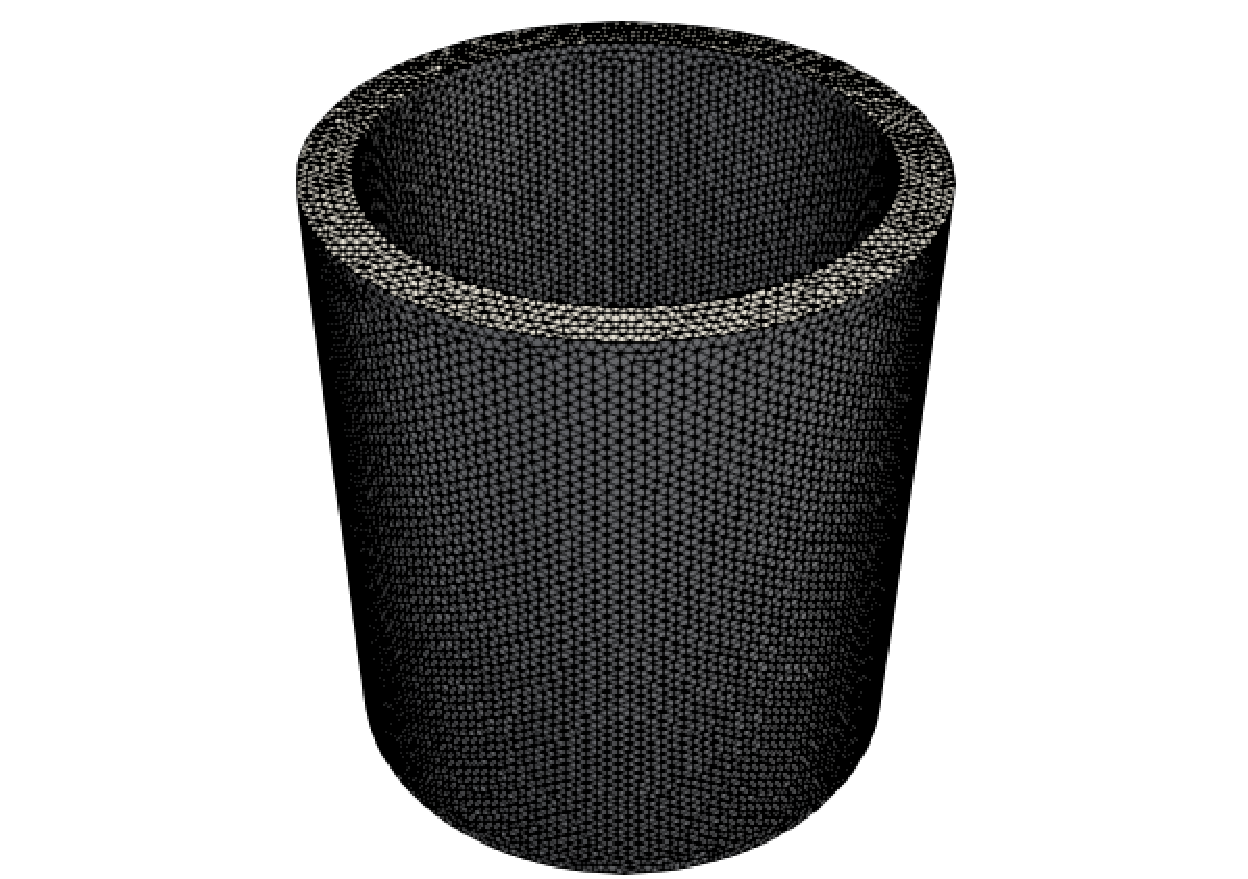}
\caption{The domain $\Omega$.}
\label{fig:first_image}
\end{minipage}\hfill\begin{minipage}{0.5\textwidth}
\centering
\includegraphics[width=\linewidth]{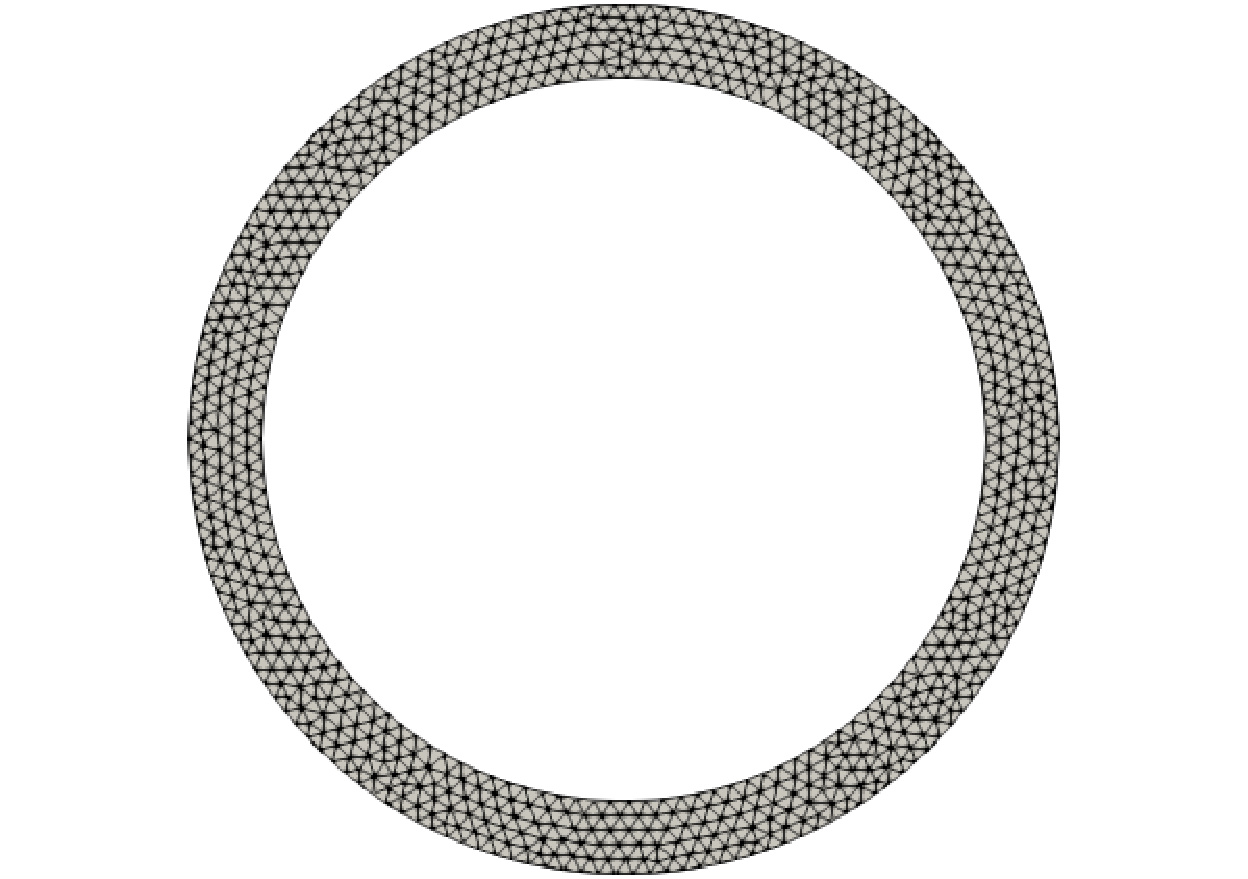}
\caption{The mesh viewed from the top.}
\label{fig:second_image}
\end{minipage}
\end{figure}

\subsubsection{3d Statistical result analysis}

We first compared the 1/2-equation velocity to both 1-equation velocities. All
3 models gave a time averaged kinetic energy (to 2 digits) of $3.6$. The L2
norms of time averaged differences $||v_{1/2eqn}-v_{1eqn}||$, $||v_{1/2eqn}%
-v_{1eqn\&l=\kappa y}||$\ were (to 2 digits) $0.17$ and $0.14$ respectively.
(These norms were calculated using nodal values in a standard way and are
known in finite element theory \cite{BS08} to be equivalent to the continuous
$L^{2}$ norms.) This yielded percent difference of respectively $4.7\%$ and
$3.9\%$:%
\[
\frac{||v_{1/2eqn}-v_{1eqn}||_{L^{2}}}{||v_{1eqn}||_{L^{2}}}\simeq0.047\text{
and }\frac{||v_{1/2eqn}-v_{1eqn\&l=\kappa y}||_{L^{2}}}{||v_{1eqn\&l=\kappa
y}||_{L^{2}}}\simeq0.039.
\]
Given the 1/2-equation model parameters were non-calibrated, these velocity
differences\ seem acceptable.

\begin{figure}[ptbh]
\centering
\includegraphics[width=0.7\linewidth]{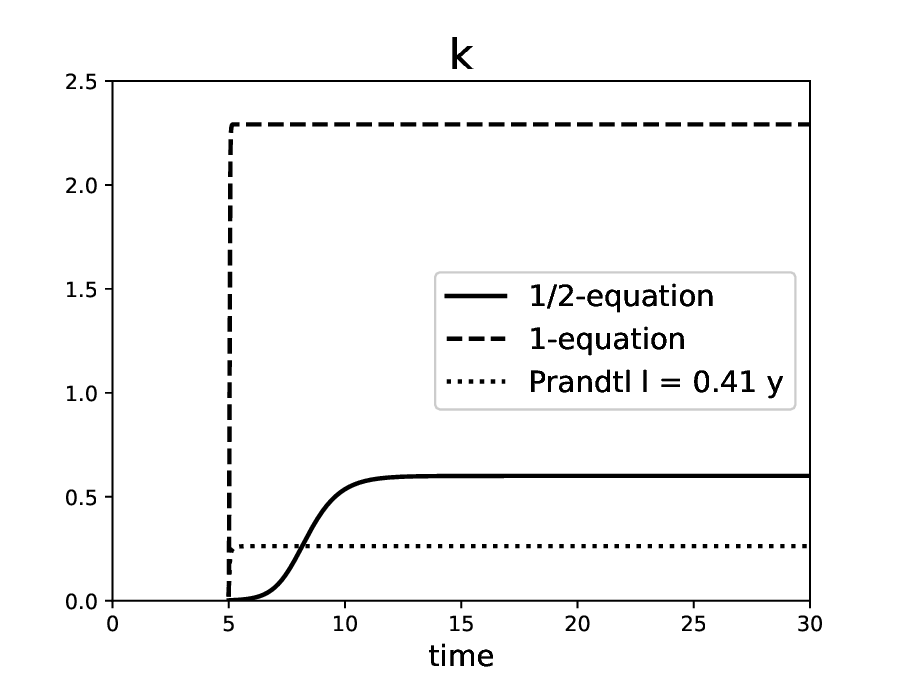}\caption{
The k value.}%
\label{fig:kn}%
\end{figure}\begin{figure}[ptb]
\centering
\includegraphics[width=0.7\linewidth]{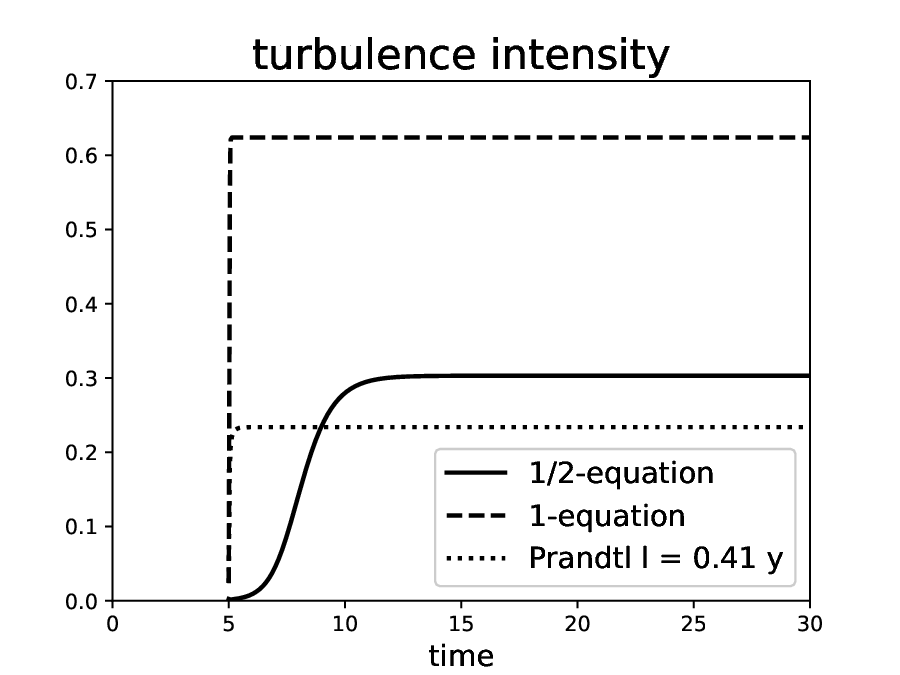}\caption{
The turbulence intensity.}%
\label{fig:turbulence_intensity}%
\end{figure}

The tests did show $k$-value differences between the non-calibrated
1/2-\-equation $k(t)$ and the space average of $k(x,t)$ for\ the 2 models.
These differences are also reflected in the computed approximations to the
turbulent intensities (as these depend on the $k$ values). In the 2d tests
below a well resolved NSE simulation was available for comparison. The 2d
results suggest that here the $k$ value of the 1-equation model is too large
due to $\nu_{T}$\ being too large in the near wall region. This suggests the
1/2-equation model results for $k(t)$ and the turbulent intensity, being
closer to the model with $l=0.41y$, are again acceptable, Figure \ref{fig:kn},
Figure \ref{fig:turbulence_intensity}. \begin{figure}[ptbh]
\centering
\includegraphics[width=0.7\linewidth]{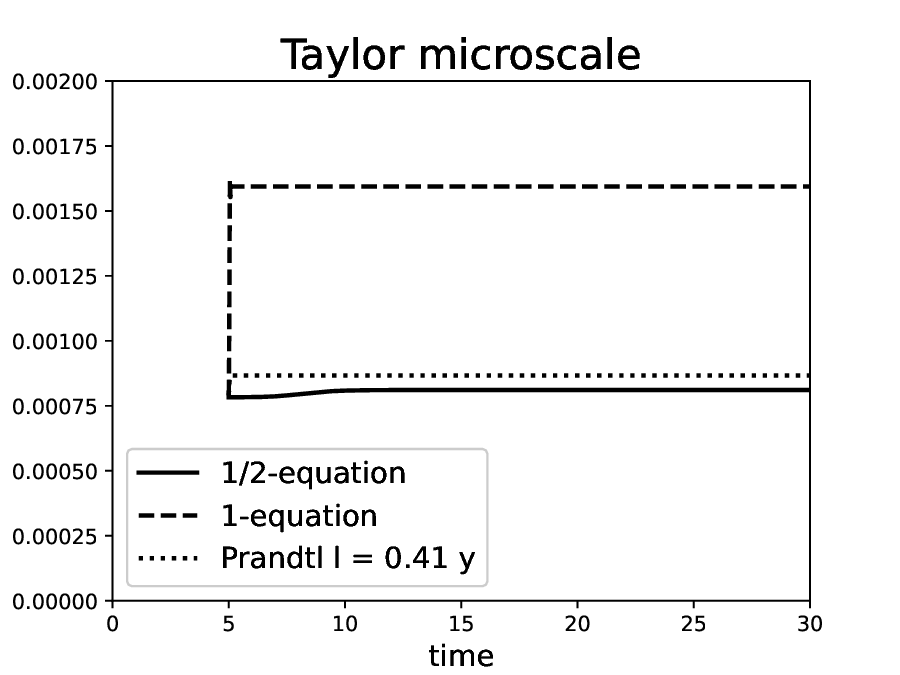}
%\end{minipage}
\caption{The Taylor microscale.}%
\label{fig:taylor_microscale}%
\end{figure}

The Taylor microscale $\lambda_{t}$ depends on velocity gradients which are
more sensitive to model parameters and mesh than velocities. Predictions of
$\lambda_{t}$ are very similar in all models from $10^{-3}$ to $1.5\times
10^{-3}$, Figure \ref{fig:taylor_microscale}. For this problem, we believe the
model with $l=0.41y$\ is more accurate than the other 1-equation model due to
its near wall asymptotics being closer to that of the Reynolds stress. Thus,
the 1/2-equation model's closeness to the former is another model success.
\begin{figure}[ptbh]
\centering
\includegraphics[width=0.7\linewidth]{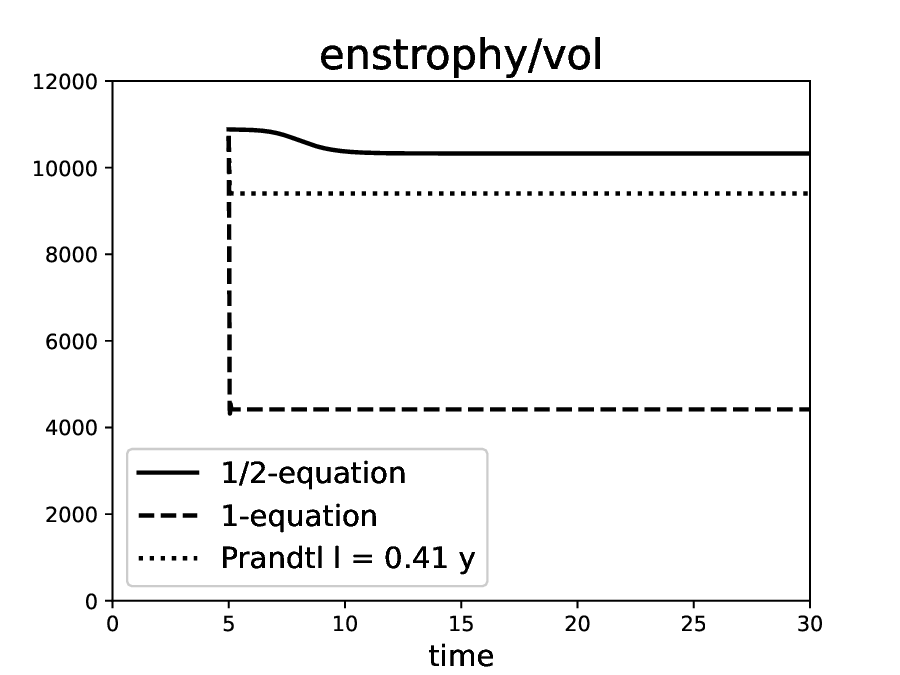}\caption{
The enstrophy over volume.}%
\label{fig:enstrophy_over_vol}%
\end{figure}

\begin{figure}[ptbh]
\begin{minipage}{0.32\textwidth}
\centering
\includegraphics[width=\linewidth]{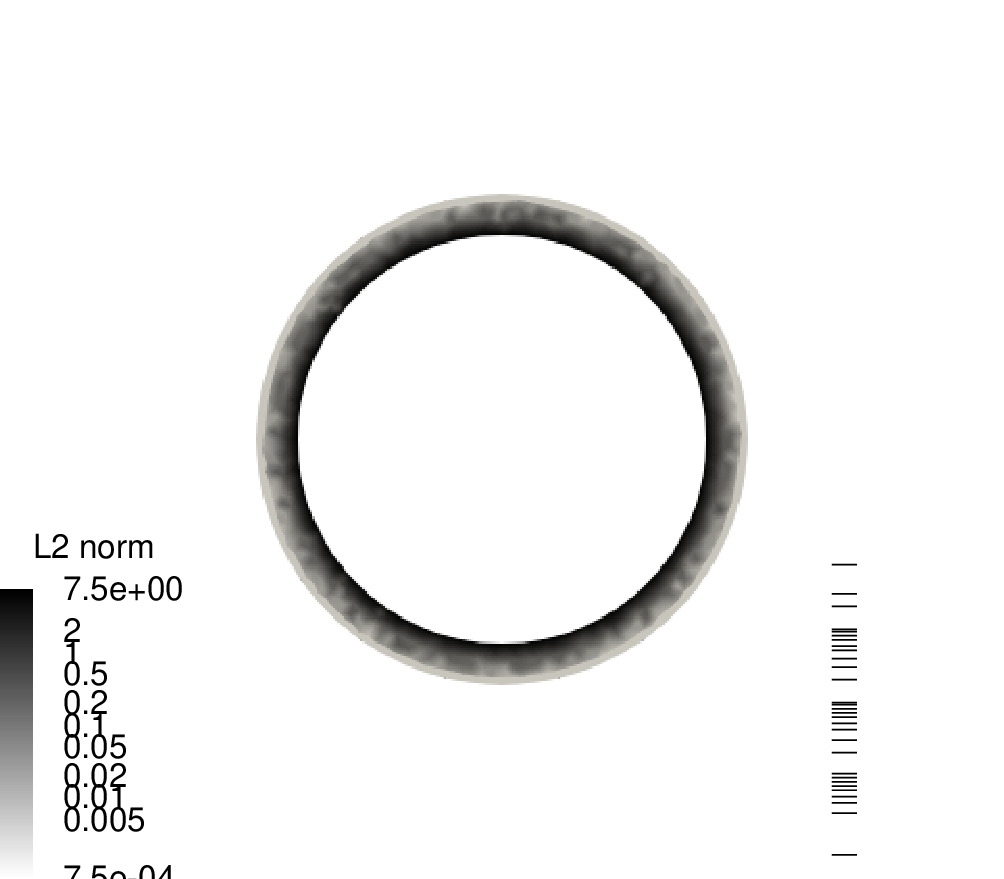}
\caption{1/2-equation:$\|u\|$ at T=30.}
\label{fig: magnitude 1/2-equation}
\end{minipage}
\begin{minipage}{0.32\textwidth}
\centering
\includegraphics[width=\linewidth]{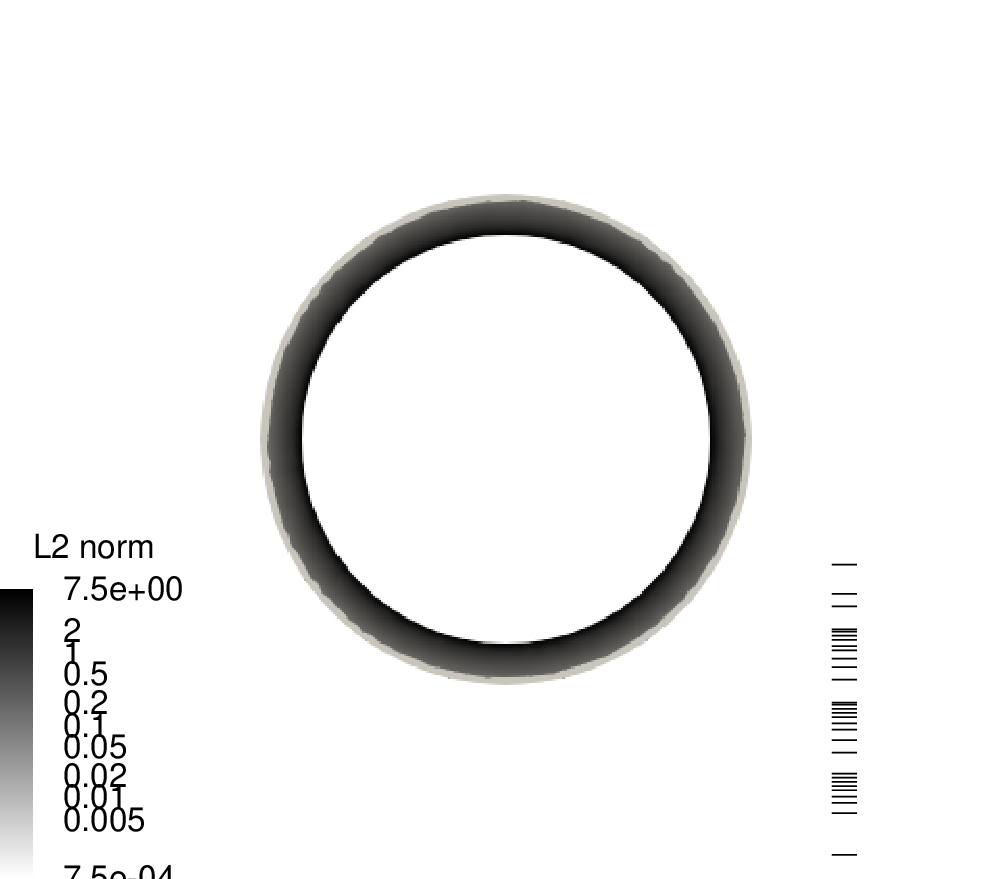}
\caption{1-equation: $\|u\|$ at T=30.}
\label{fig: magnitude 1-equation}
\end{minipage}
\begin{minipage}{0.3\textwidth}
\centering
\includegraphics[width=\linewidth]{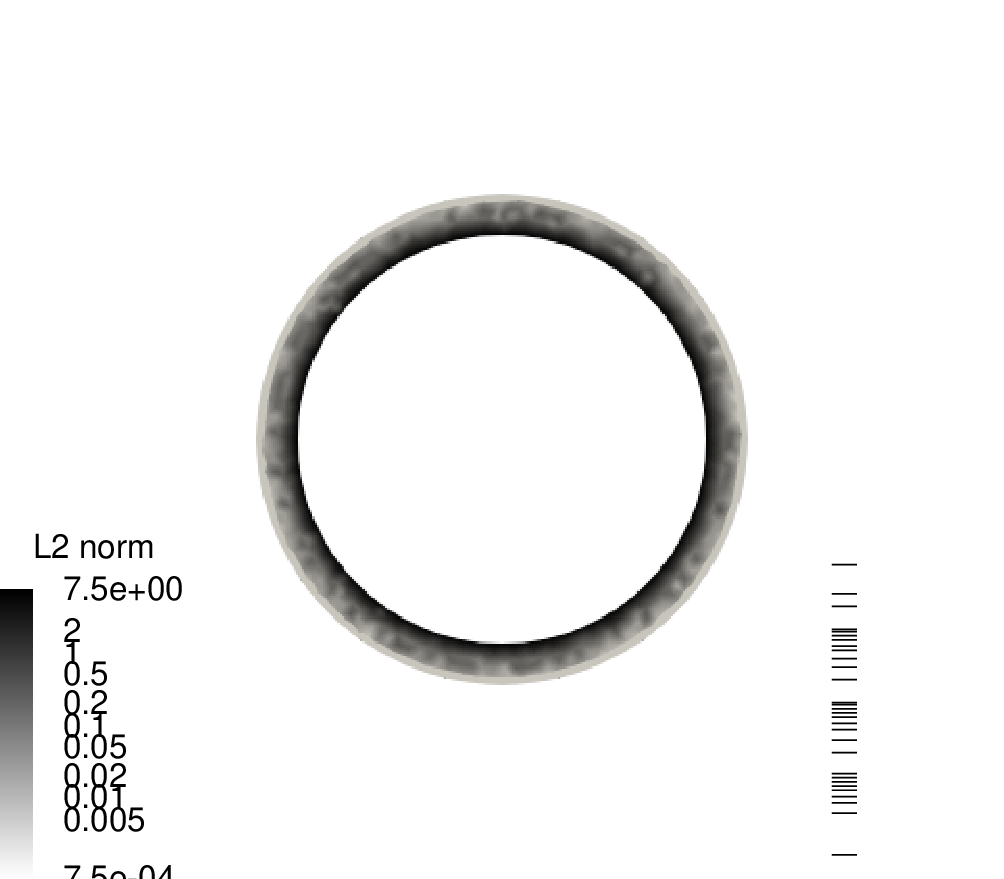}
\caption{Prandtl: $\|u\|$ at T=30.}
\label{fig: magnitude prandtl}
\end{minipage}\hfill\end{figure}

\begin{figure}[ptbh]
\begin{minipage}{0.32\textwidth}
\centering
\includegraphics[width=\linewidth]{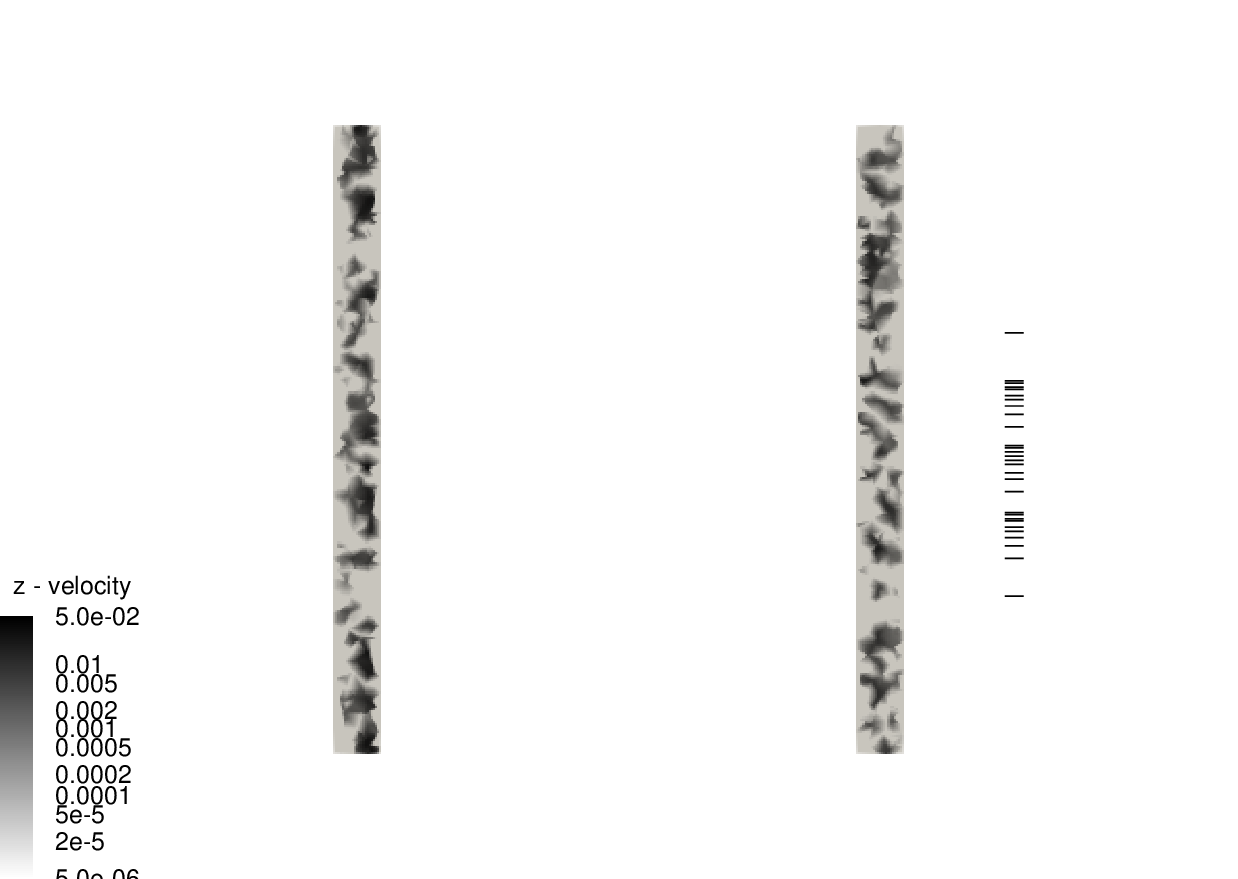}
\caption{1/2-equation: z component of the velocity at T=30.}
\label{fig: z 1/2-equation}
\end{minipage}
%\end{figure}
%\begin{figure}[htbp]
\begin{minipage}{0.32\textwidth}
\centering
\includegraphics[width=\linewidth]{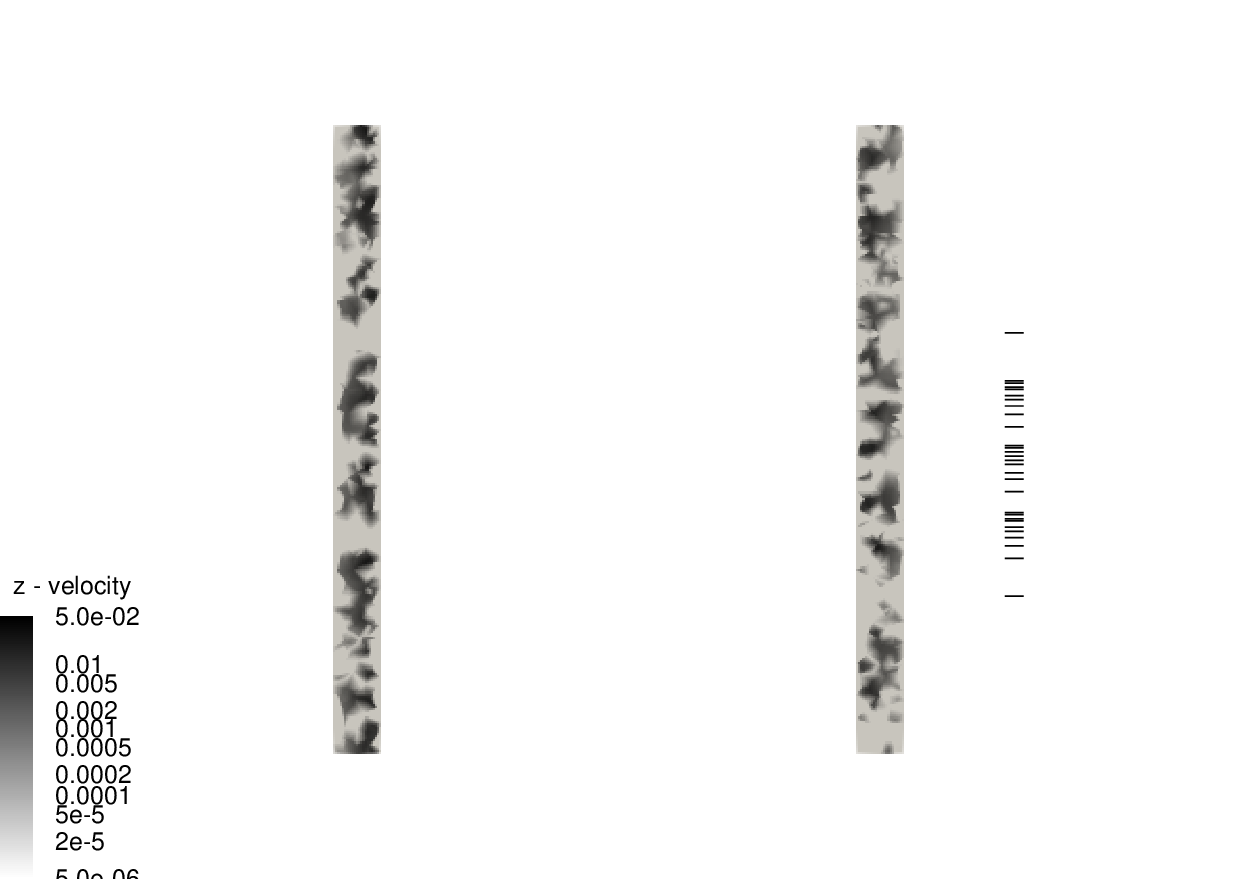}
\caption{1-equation: z component of the velocity at T=30.}
\label{fig: z 1-equation}
\end{minipage}
\begin{minipage}{0.32\textwidth}
\centering
\includegraphics[width=\linewidth]{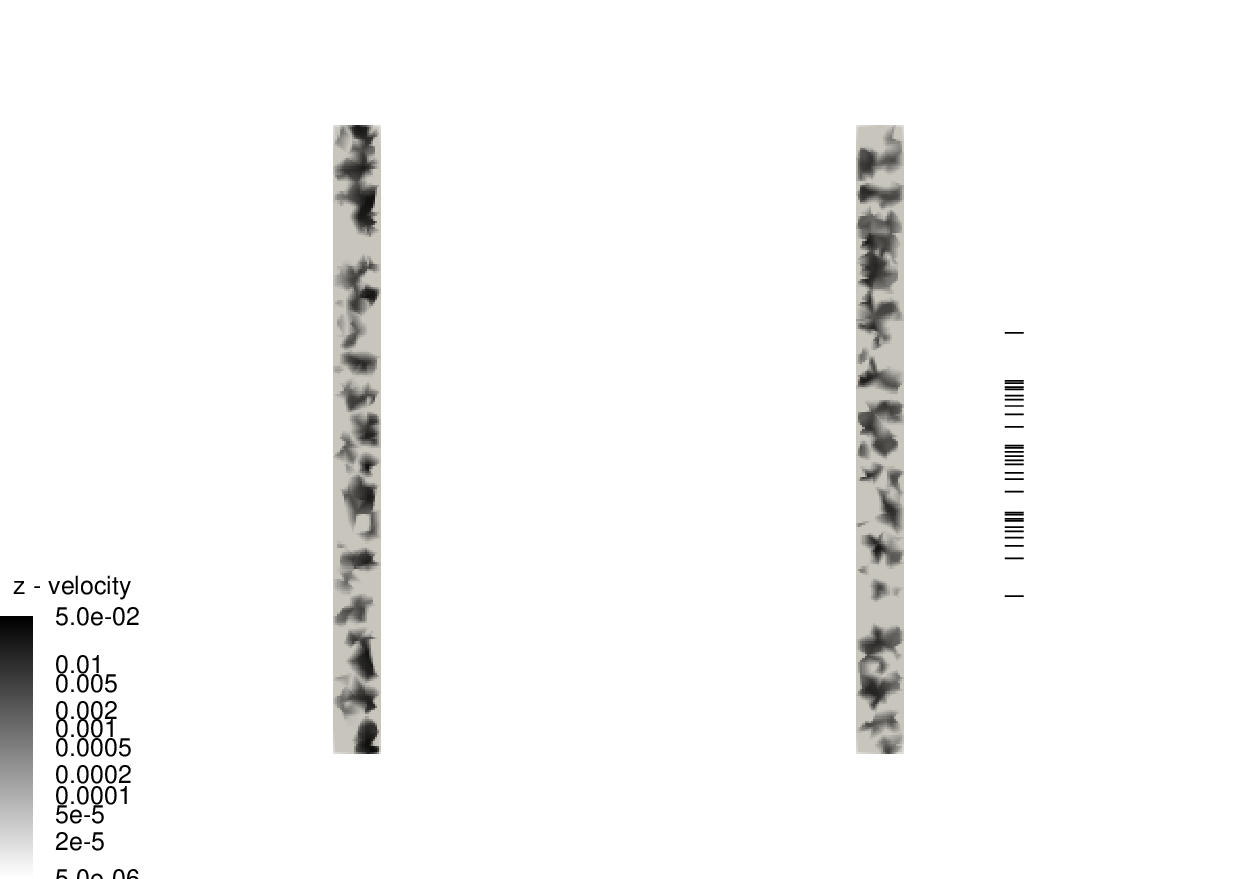}
\caption{Prandtl: z component of the velocity at T=30.}
\label{fig: z prandtl}
\end{minipage}\hfill
\end{figure}

The enstrophy values indicate significant rotational motions. The previous
results suggest the $k(x,t)$ values for the 1-equation model with $l=\sqrt
{2k}\tau$ are too large for this problem. This makes $\nu_{T}$ too large and
the model velocity over diffused. Thus lower enstrophy is expected.\ In Figure
\ref{fig:enstrophy_over_vol}, 1/2-equation model has enstrophy close to the
model solution with $l=\kappa y$ and above the 1-equation model with
$l=\sqrt{2k}\tau$ . The magnitude of velocity at time $T=30$ for all the
models is presented in Figure \ref{fig: magnitude 1/2-equation},
\ref{fig: magnitude 1-equation} and \ref{fig: magnitude prandtl}. The $z$
components of models' velocity are plotted in Figure \ref{fig: z 1/2-equation}%
, \ref{fig: z 1-equation}, \ref{fig: z prandtl}. In these we can see vertical
rotations consistent with irregular Taylor cells. 

\subsection{The 2d test problem}

For the 2d tests of flow between offset circles we selected $\kappa=0.41$ in
the multiplier $(\kappa y/L)^{2}$. Since this problem is 2d, we were able to
perform a well resolved NSE simulation for comparison. We compared the
1/2-equation model velocity statistics to the 1-equation model statistics with
$\nu_{T}=\sqrt{2}\mu k(x,t)\tau$ and with velocity statistics computed from
the well-resolved NSE simulation. The other details of the 2d tests are as
follows. The computational domain is a disk with a smaller off center obstacle
inside.
\[
\Omega=\{(x,y):x^{2}+y^{2}\leq r_{1}^{2}\cap(x-c_{1})^{2}+(y-c_{2})^{2}\geq
r_{2}^{2}\}
\]
where we set $r_{1}=1,\quad r_{2}=0.1,\quad c=(c_{1},c_{2})=(\frac{1}{2},0)$.
No-slip boundary conditions are imposed on both circles. The flow is driven by
a counterclockwise force $f(x,y,t)=(-4y\min\left(  t,1\right)  (1-x^{2}%
-y^{2}),4x\min\left(  t,1\right)  (1-x^{2}-y^{2}))$. We set $\tau=0.1$,
$\mu=0.55$, $\nu=10^{-4}$, $L=1$, $U_{\max}=1$ and $\text{Re}=\frac{UL}{\nu}$.
The final time is $T=15$. The k-equation is initialized at $t^{\ast}=1$.

\textbf{Initial and boundary conditions:} For the 1-equation model and the
1/2-equation model, we choose initialization for the 2 $k-$equations as in
\cite{KLS22}: $t^{\ast}=1$\ and
\[
k\left(x,1\right)  =\frac{1}{2\tau^{2}}l^{2}\left(  x\right),
l\left(  x\right)  =\min\left\{  \kappa y,0.082Re^{-1/2}\right\} \text{ \&
}\]
\[k(1)=\frac{1}{|\Omega|}\frac{1}{2\tau^{2}}\int_{\Omega}l\left(  x\right)
^{2}dx.
\]
The boundary condition for the k-equation is homogeneous Dirichlet.

\textbf{Discretization:} We employ the Taylor-Hood $\left(  P2-P1\right)  $
finite element pair for approximating the velocity and pressure and $P1$
Lagrange element for the TKE equation. We choose the timestep $\Delta t=0.01$
and use the backward Euler time discretization. The mesh is generated by the
Delaunay triangular method with 40 mesh points on the outer circle and 20 mesh
points on the inner circle. This mesh has the longest edge $\max_{e}%
h_{e}=0.208201$ and the shortest edge $\min_{e}h_{e}=0.0255759$. For the well
resolved NSE solve, we use a finer mesh with 80 mesh points on the outer
circle and 60 mesh points on the inner circle, extended by a Delaunay
triangulation. This mesh has the longest edge $\max_{e}h_{e}=0.108046$ and the
shortest one $\min_{e}h_{e}=0.0110964$. The 2d tests were performed with
FreeFEM++, Hecht \cite{H12}.

\subsubsection{2d Statistical result analysis}

The space average of the 1-equation model's $k(x,t)$ was larger than the
1/2-equation model's $k(t)$ as in 3d. This is likely because of the difference
between the sizes of the two models $\nu_{T}$\ values near the inner
disk\footnote{Options to correct the 1-equation model include near wall
clipping \cite{KLS21} or rescaling \cite{KLS22} or damping functions
\cite{P17}. These were not done because we test here the 1/2 equation model.}.
This led to the question of which model's velocity statistics were more
accurate. For this reason we performed the well resolved NSE simulation.

\begin{figure}[ptbh]
\centering
\includegraphics[width=0.7\linewidth]{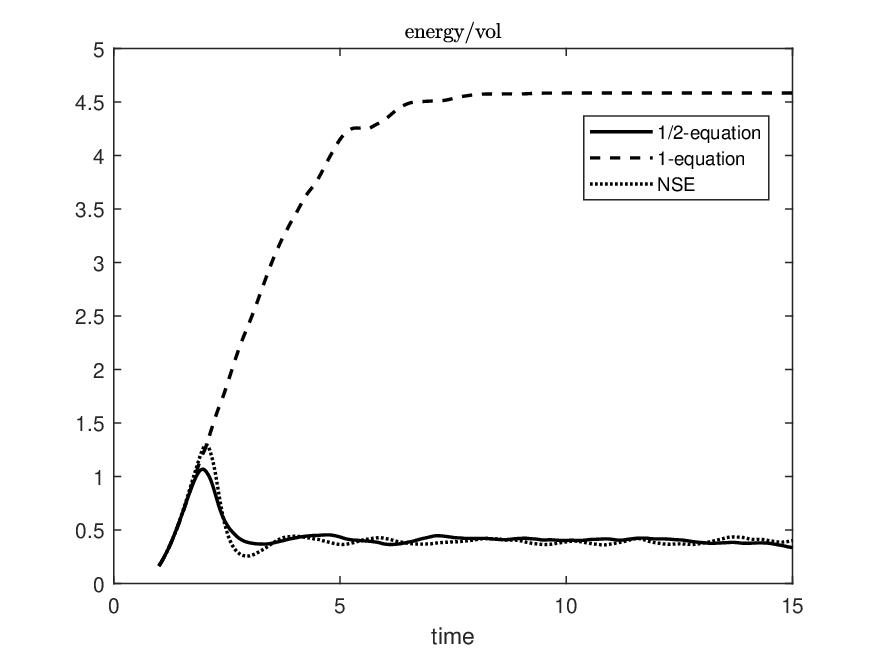}\caption{
The kinetic energy over volume.}%
\label{fig: 2d_KineticEnergy}%
\end{figure}

\begin{figure}[ptbh]
\centering
\includegraphics[width=0.7\linewidth]{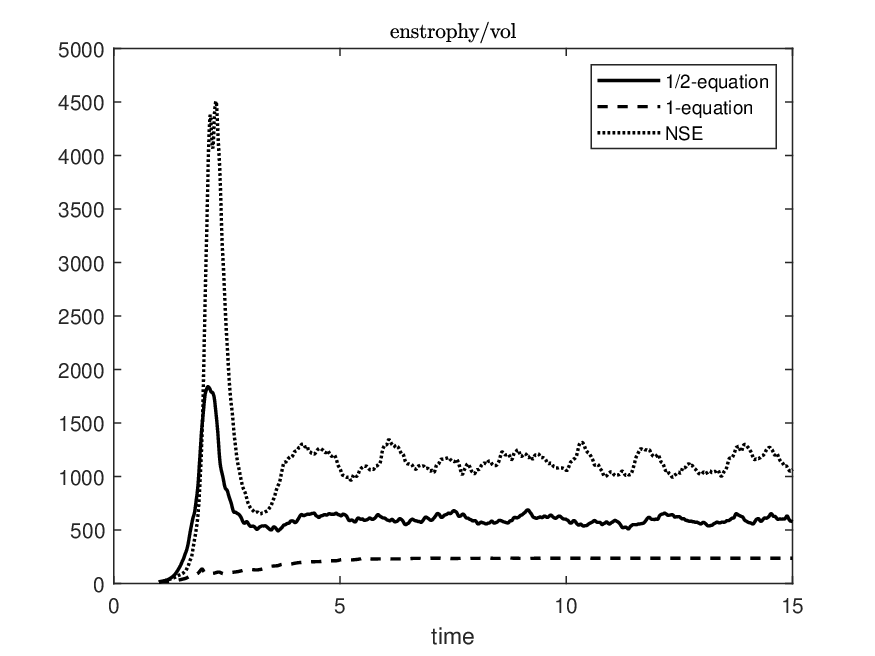}\caption{
The enstrophy over volume.}%
\label{fig: 2d_Enstrophy}%
\end{figure}

\begin{figure}[ptbh]
\centering
\includegraphics[width=0.7\linewidth]{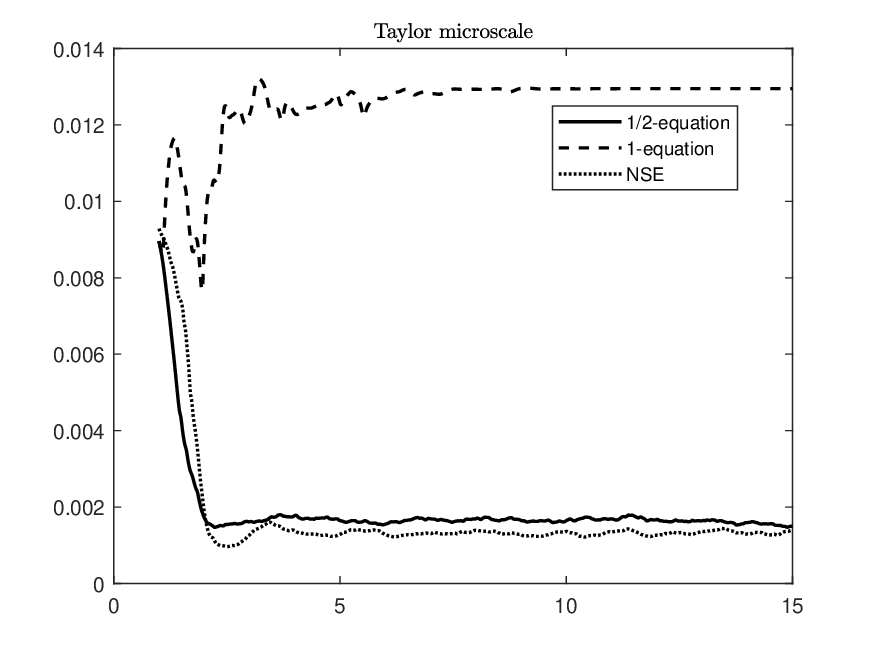}\caption{
The Taylor microscale.}%
\label{fig: 2d_TaylorM}%
\end{figure}

Fig \ref{fig: 2d_KineticEnergy}, \ref{fig: 2d_Enstrophy}, and
\ref{fig: 2d_TaylorM} present the comparison of the evolution of the
respective kinetic energies, enstrophy and Taylor microscales. We observe that
the kinetic energy of the 1/2-equation model is slightly less than that of the
well resolved NSE test but closely tracks\ it's behavior. The 1-equation
model's kinetic energy is clearly incorrect. The same behavior was observed
for the enstrophy and Taylor microscales in Figures \ref{fig: 2d_Enstrophy},
and \ref{fig: 2d_TaylorM}. 

\section{Conclusions}

Due to the computational costs of DNS and LES, RANS and URANS models are still
widely used. This suggests two fruitful directions of URANS research: lowering
simulation costs preserving current accuracy and raising accuracy at current
simulation costs. The 1/2-equation model herein aims at the former. The model
derivation, analysis and tests indicate the 1/2-equation model is worthy of
further study and the idea behind it of further development. The 1/2-equation
model (\ref{eqHalfEqnModel}) produced velocity-statistics comparable to the
same velocity-statistics for 1-equation models in our tests. No model is
perfect so further tests delineating failure modes would be useful. The next
interesting tests include flows with time-varying body forces and with
interior shear layers. There are also many parallel analytical questions.

Our longer term motivation was to use a similar idea to simplify more complex
models such as 2-equation models. In these the TKE equation is well grounded
in mechanics but the second equation, used to indirectly determine the
turbulence length scale, is often a product of optimism, data fitting and
experience-informed intuition. For these an simplified model for $l(t)$ (thus
a 3/2-equation model) is an interesting possibility to explore, building on
work here in the most basic case. 

\begin{section}{Acknowledgement}
This research herein of W. Layton and Rui Fang was supported in part by the
NSF under grant DMS\ 2110379. WE also gratefully acknowledge support of the
University of Pittsburgh Center for Research Computing through the resources
provided on the SMP cluster. The author Weiwei Han was partially supported by
the Innovative Leading Talents Scholarship established by Xi'an Jiaotong University.
\end{section}





\end{document}